\documentclass[a4paper,12pt]{article}
 \usepackage[utf8,latin1]{inputenc} 
\usepackage[T1]{fontenc} 
\usepackage{mathtools,amsmath}
\usepackage{amssymb}
\usepackage{array}

\begin{document}

\title{\bf A result related to the Sendov conjecture}
\author{{\bf Robert DALMASSO}\\\\ Le Galion - B,\\ 
33 Boulevard Stalingrad,\\
 06300 Nice, France}

\date{}
\maketitle
\footnote[0]{\sl

E-mail address: {\rm robert.dalmasso51@laposte.net}


}

\allowdisplaybreaks[4]

\bigskip

\hrule

\bigskip

\noindent{\bf Abstract.} The Sendov conjecture asserts that if $p(z) = \prod_{j=1}^{N}(z-z_j)$ is a polynomial with zeros $|z_j| \leq 1$, then each disk $|z-z_j| \leq 1$ contains a zero of $p'$. Our purpose is the following: Given a zero $z_j$ of order $n \geq 2$, determine whether there exists $\zeta \not= z_j$ such that $p'(\zeta) = 0$ and $|z_j - \zeta| \leq 1$. In this paper we present some partial results on the problem.

\bigskip

\hrule

\bigskip

The well-known Sendov conjecture (\cite{h} Problem 4.5) asserts that if $p(z) = \prod_{j=1}^{n}(z - z_j)$ is a polynomial with $|z_j| \leq 1$ ($1 \leq j \leq n$), then each disk $|z - z_j|\leq 1$ ($1 \leq j \leq n$) contains a zero of $p'$. Notice that by the Gauss-Lucas theorem the zeros $w_k$ ($1 \leq k \leq n-1$) of $p'$ lie in the closed convex hull of the zeros of $p$, hence $|w_k| \leq 1$ for $1 \leq k \leq n-1$. 

This conjecture has been verified for polynomials of degree $n \leq 8$ or for arbitrary degree $n$ if there are at most eight distinct roots: See Brown and Xiang \cite{bx} and the references therein. It is also true in general ($n \geq 2$)  when $p(0) = 0$ (\cite{sc}). The Sendov conjecture is true with respect to the root $z_j$ of $p$ if $|z_j| = 1$ (\cite{ru}). An improved version of this result was later given (\cite{grr}). Recently (\cite{de}) it has been verified when $n$ is larger than a fixed integer depending on the root $z_j$ of $p$ (see also \cite{ch}). Finally Tao \cite{t} proved the conjecture for sufficiently high degree polynomials. We refer the reader to M. Marden \cite{m} and Bl. Sendov \cite{s} for further information and bibliographies.

\bigskip

In this paper we are interested in the following problem.

\bigskip

\noindent{\bf Problem:} Let
\begin{equation}
p(z) = (z-a)^n\prod_{j=1}^{k}(z-z_j)^{n_j} \,\, , \,\, z \in \mathbb{C}\, ,
\label{eq:eq1}
\end{equation}
where $|a| \leq 1$, $z_j \not= a$ and $|z_j| \leq 1$ for $j = 1,\cdots ,k$ ($k \geq 1$) with $z_i \not= z_j$ for $i \not= j$ and $n$, $n_j \geq 1$ for $j = 1,\cdots ,k$.

Suppose that $n \geq 2$. We want to determine whether there exists $\zeta \not= a$ such that $p'(\zeta) = 0$ and $|a - \zeta| \leq 1$.

\bigskip

Clearly we can assume that $0 \leq a \leq 1$. The answer is obviously positive if $a = 0$ and also  if there exists $j \in \{1,\cdots,k\}$ such that $n_j \geq 2$ and $|a-z_j| \leq 1$. 

However the answer to the problem is not positive in general. Indeed we have the following lemma.

\bigskip

\noindent{\bf Lemma.} {\sl Let $k$ and $r$ be positive integers and let $|a| \leq 1$, $|b| \leq 1$ be such that $|a-b| > 1$. Define
\begin{displaymath}
p(z) = (z-a)^k(z-b)^r \,\, , \,\, z \in \mathbb{C}\, .
\end{displaymath}
Then there exists $\zeta \not= a$ such that $p'(\zeta) = 0$ and $|a - \zeta| \leq 1$ if and only if $k \leq r$.}

\bigskip

\noindent{\sl Proof.} We have

\begin{displaymath}
p'(z) = (k+r)(z-a)^{k-1}(z-b)^{r-1}\Big (z-\displaystyle\frac{ra+kb}{k+r}\Big) ,
\end{displaymath}
and
 \begin{displaymath}
\Big |\displaystyle\frac{ra+kb}{k+r} - a \Big| = \displaystyle\frac{k}{k+r}|a-b| > 1\, ,
\end{displaymath}
if
 \begin{equation}
k(|a-b|-1) > r \, .
\label{eq:eq2}
\end{equation}
\eqref{eq:eq2} implies that $k > r$ since $|a-b| \leq 2$. On the other hand, if $2 \leq k \leq r$ $\zeta = (ra+kb)/(k+r) \not= a$ is such that $p'(\zeta) = 0$ and $|a - \zeta| \leq 1$ since $k/(k+r) \leq 1/2$ and $|a-b| \leq 2$. The proof is complete.

\bigskip

We have
\begin{displaymath}
p'(z) =  m(z-a)^{n-1}\prod_{j=1}^{k}(z-z_j)^{n_j -1}\prod_{j=1}^{k}(z-w_j) \,\, , \,\, z \in \mathbb{C}\, ,
\end{displaymath}
where
\begin{displaymath}
m = n + \sum_{j=1}^{k}n_j\, .
\end{displaymath}

\bigskip

Now we can state our  first result.

\bigskip

\noindent{\bf Theorem 1.} {\sl Let $p$ be as in \eqref{eq:eq1} with $ a = 1$, $z_j \not= 1$ and $|z_j| \leq 1$ for $j = 1,\cdots,k$ ($k \geq 1$) with $z_i \not= z_j$ for $i \not= j$. Let $n \geq 1$ and $n_j \geq 1$ for $j = 1,\cdots,k$. Then there exists $\zeta \not= a$ such that $p'(\zeta) = 0$ and 
\begin{displaymath}
\displaystyle \Big |\zeta - \frac{m-n}{m + (k-1)n} \Big| \leq \frac{kn}{m + (k-1)n}  \,\, .
\end{displaymath}}

\bigskip

\noindent{\sl Proof.} Let
\begin{displaymath}
p(z) = (z-1)^{n}q(z)\, .
\end{displaymath}
We have
\begin{displaymath}
p'(z) = (z-1)^{n-1}(nq(z) + (z-1)q'(z))  \,\, , \,\, z \in \mathbb{C}\, ,
\end{displaymath}
and
\begin{displaymath}
p''(z) = (z-1)^{n-2}(n(n-1)q(z) + 2n(z-1)q'(z) + (z-1)^{2}q''(z)) \,\, , \,\, z \in \mathbb{C}\, .
\end{displaymath}
Since $q(1) \not= 0$ there exists $\varepsilon_0 > 0$ such that $nq(1-\varepsilon) -\varepsilon q'(1-\varepsilon) \not= 0$ for $\varepsilon \in (0,\varepsilon_0)$. Then, for $\varepsilon \in (0,\varepsilon_0)$, we can write
\begin{displaymath}
\begin{array}{lcl}
 \displaystyle  \frac{p''(1- \varepsilon)}{p'1- \varepsilon)} & = &  \displaystyle  \frac{-n(n-1)q(1- \varepsilon) + 2n\varepsilon q'(1- \varepsilon) - \varepsilon^2 q''(1- \varepsilon)}{n\varepsilon q(1- \varepsilon) - \varepsilon^2 q'(1- \varepsilon)} \\ \\& = & \displaystyle - \frac{n-1}{\varepsilon} + \frac{n+1}{n}\frac{q'(1- \varepsilon)}{q(1- \varepsilon)} + A(\varepsilon) \, \, ,
 \end{array}
\end{displaymath}
where $A(\varepsilon) \to 0$ when $\varepsilon \to 0$. On the other hand we have
\begin{displaymath}
 \displaystyle  \frac{p''(1- \varepsilon)}{p'1- \varepsilon)} = \displaystyle - \frac{n-1}{\varepsilon} + \sum_{j=1}^{k}\frac{n_j -1}{ 1 - \varepsilon - z_j} + \sum_{j=1}^{k}\frac{1}{ 1 - \varepsilon - w_j} \, \, .
\end{displaymath}
We deduce that
\begin{displaymath}
\begin{array}{lcl}
 \displaystyle \sum_{j=1}^{k}\frac{n_j -1}{ 1 - \varepsilon - z_j} + \sum_{j=1}^{k}\frac{1}{ 1 - \varepsilon - w_j} & = &  \displaystyle  \frac{n+1}{n}\frac{q'(1-\varepsilon)}{q(1- \varepsilon)} + A(\varepsilon) \\ \\& = & \displaystyle \frac{n+1}{n}\sum_{j=1}^{k}\frac{n_j}{ 1 - \varepsilon - z_j}  + A(\varepsilon) \, \, .
 \end{array}
\end{displaymath}
Letting $\varepsilon \to 0$ we get
\begin{displaymath}
 \displaystyle   \sum_{j=1}^{k}\frac{1}{ 1 - w_j} = \frac{1}{n}\sum_{j=1}^{k}\frac{n_j +n}{ 1 - z_j}   \, \, .
\end{displaymath}
Since ${\rm Re} \, 1/(1-z_j) \geq 1/2$ we obtain
\begin{displaymath}
 \displaystyle {\rm Re} \, \sum_{j=1}^{k}\frac{1}{ 1 - w_j}  \geq \frac{m + (k-1)n}{2n}\, \, .
\end{displaymath}
It follows that, for at least one $j$, we have
\begin{equation}
 \displaystyle {\rm Re} \, \frac{1}{ 1 - w_j}  \geq \frac{m + (k-1)n}{2kn} \, \, .
 \label{eq:eq3}
 \end{equation}
This is equivalent to
\begin{displaymath}
\Big |w_j -\frac{m-n}{m + (k-1)n}\Big| \leq \frac{kn}{m + (k-1)n}  \,\, ,
\end{displaymath}
as required.

\bigskip

\noindent{\bf Corollary 1.} {\sl Let $p$ be as in \eqref{eq:eq1} with $ a = 1$, $z_j \not= 1$ and $|z_j| \leq 1$ for $j = 1,\cdots,k$ ($k \geq 1$) with $z_i \not= z_j$ for $i \not= j$. Let 
\begin{displaymath}
n \geq 2 \quad  {\textrm and} \quad \sum_{j=1}^{k}n_j \geq kn \, \, .
\end{displaymath}
Then there exists $\zeta \not= a$ such that $p'(\zeta) = 0$ and $|\zeta - 1/2| \leq 1/2$.}

\bigskip

\noindent{\sl Proof.} As in the proof of Theorem 1 we obtain \eqref{eq:eq3}. Since $m + (k-1)n \geq 2kn$ we have 
\begin{displaymath}
 \displaystyle {\rm Re} \, \frac{1}{ 1 - w_j}  \geq 1 \, \, .
 \end{displaymath}
This is equivalent to
\begin{displaymath}
|w_j -\frac{1}{2}| \leq \frac{1}{2}  \,\, .
\end{displaymath}

\bigskip

In the next theorem we give another  partial result. Let $\varphi$ be defined by
\begin{displaymath}
\displaystyle \varphi(a) = (1-a)^k - \frac{a(m-n)+n}{m} \quad , \, \, 0 \leq a \leq 1 \, \, . 
\end{displaymath}
There exists $a_0 \in (0,1)$ such that $\varphi(a_0) = 0$ and $\varphi(a) > 0$ for $a \in [0,a_0)$.

\bigskip

\noindent{\bf Theorem 2.} {\sl Let $p$ be as in \eqref{eq:eq1} with $0 < a \leq a_0$, $z_j \not= a$  for $j = 1,\cdots,k$ ($k \geq 1$) with $z_i \not= z_j$ for $i \not= j$. Then there exists $\zeta \not= a$ such that $p'(\zeta) = 0$ and $|a - \zeta| \leq 1$.}

\bigskip

\noindent{\sl Proof.} We can write
\begin{displaymath}
p'(z) =  q(z)(z-a)^{n-1}\prod_{j=1}^{k}(z-z_j)^{n_j -1} \,\, , \,\, z \in \mathbb{C}\, ,
\end{displaymath}
where 
\begin{displaymath}
q(z) =  n\prod_{j=1}^{k}(z-z_j) + (z-a)\sum_{j=1}^{k}n_j\prod_{\substack{i=1 \\ i\neq j}}^{k}(z-z_i)  = m\prod_{j=1}^{k}(z-w_j)\,\, , \,\, z \in \mathbb{C}\, .
\end{displaymath}
We have
\begin{displaymath}
|q(0)| =  m\prod_{j=1}^{k}|w_j| \leq n + a\sum_{j=1}^{k}n_j = a(m-n)+n\, .
\end{displaymath}
Therefore there exists $j \in \{1,\cdots,k\}$ such that
\begin{displaymath}
\displaystyle |w_j| \leq  \Big (\frac{a(m-n)+n)}{m}\Big)^{\frac{1}{k}}\, .
\end{displaymath}
Then for $a \in [0,a_0]$ we get
\begin{displaymath}
\displaystyle |a-w_j| \leq a+|w_j| \leq a + \Big(\frac{a(m-n)+n)}{m}\Big)^{\frac{1}{k}} \leq 1\, .
\end{displaymath}
The proof is complete.
\bigskip

\noindent {\sl Remark 1}. If $n_1 + \cdots + n_k \geq kn$, we have 
\begin{displaymath}
 \displaystyle \varphi(a) \geq (1-a)^k -  \, \frac{ak+1}{k+1}  \, \, ,
 \end{displaymath}
and we easily show that $\varphi(1/k) > 0$ when $k \geq 7$, hence $a_0 > 1/k$ for $k \geq 7$.

\bigskip

Our last result is the following theorem.

\bigskip

\noindent{\bf Theorem 3 .} {\sl Let $p$ be as in (1) with $0 < a \leq 1$, $z_j \not= a$, $|z_j| \leq 1$ for $j = 1,\cdots ,k$ ($k \geq 1$) and  $z_i \not= z_j$ for $i \not= j$. Let $n \geq 2$ and $n_j \geq 1$ for $j = 1,\cdots,k$. Assume the following hypotheses:

(i) $n_j = r_jn + s_j(n)$  where $0 \leq s_j(n) \leq n -1$ for $j = 1,\cdots,k$ and $r_1, \cdots, r_k$ do not depend on $n$;

(ii) $r_j + s_j(n) \geq 1$ for $j = 1,\cdots,k$ and $r_j \geq 1$ for some $j \in \{1,\cdots,k\}$;

(iii) $s_j(n)/n \to 0$ as $n \to \infty$ for $j \in \{1,\cdots,k\}$;

(iv) for
\begin{displaymath}
X = \{j \in \{1,\cdots, k\}; r_j \geq 1\} \, , \quad Y = \{j \in \{1,\cdots, k\}; s_j(n) \geq 1\}\, ,
\end{displaymath}
and
\begin{displaymath}
P(z) = (z-a)\prod_{j\in X}(z-z_j)^{r_j} \,\, , \,\, z \in \mathbb{C}\, ,
\end{displaymath}
there exists $\zeta \in \overline{D(0,1)}$  such that $P'(\zeta) = 0$ and $|a-\zeta| < 1$.

Then there exists an integer $n_0$ such that for $n \geq n_0$ there exists $\zeta \not= a$ such that $p'(\zeta) = 0$ and $|a-\zeta| < 1$.}

\bigskip

\noindent {\sl Proof.}
We can assume that $|a-z_j| \geq 1$ for $j = 1,\cdots,k$. 

1)  If $s_j(n) = 0$ for $j = 1,\cdots,k$, then $X = \{1,\cdots,k\}$ by (ii) and $p = P^{n}$ where 
\begin{displaymath}
P(z) =  (z-a)\prod_{j =1}^{k}(z-z_j)^{r_j} \,\, , \,\, z \in \mathbb{C}\, .
\end{displaymath}
Then the result follows from (iv) with $n_0 = 2$.

2) If there exists $j \in \{1,\cdots,k\}$ such that $s_j(n) \geq 1$ we can write $p = P^{n}q_n$ where 
 \begin{displaymath}
P(z) =  (z-a)\prod_{j\in X}(z-z_j)^{r_j} \,\, , \,\, z \in \mathbb{C}\, ,
\end{displaymath}
 and
 \begin{displaymath}
q_n(z) = \prod_{j\in Y}(z-z_j)^{s_j(n)} \,\, , \,\, z \in \mathbb{C}\, .
\end{displaymath}
We have
\begin{displaymath}
p' = P^{n-1}(nP'q_n + Pq_n')\, ,
\end{displaymath}
where $nP'q_n + Pq_n'$ is a polynomial of degree $r_1 + \cdots + r_k + s_1(n)+\cdots +s_k(n)$. Denote by $\zeta_i$,  $i \in \{1, \cdots,|X|\}$ the zeros of $P'$  such that $\zeta_i \notin \{z_j;\, j \in X\}$ for $i \in \{1, \cdots,|X|\}$ . We have 
\begin{displaymath}
\begin{array}{lcl}
 nP'q_n + Pq_n' & = & \displaystyle \prod_{j=1}^{k}(z-z_j)^{r_j+s_j(n) - 1}\Big(n\Big(1 + \sum_{j=1}^{k}r_j \Big)\prod_{j = 1}^{|X|}(z-\zeta_j)\\ \\
 & & \times\displaystyle\prod_{j\in Y\backslash X}(z-z_j)\\ \\
 &  & + \displaystyle (z-a)\Big (\sum_{j=1}^{k}s_j(n)\prod_{i \in Y\backslash\{j\}}(z-z_i)\Big)\prod_{j\in X\backslash Y}(z-z_j)\Big)\, \, .
 \end{array}
\end{displaymath}
Let
\begin{displaymath}
\begin{array}{lcl}
f_n(z) & = & \displaystyle \prod_{j = 1}^{|X|}(z-\zeta_j)\prod_{j\in Y\backslash X}(z-z_j) 
 + \displaystyle \frac{z-a}{n(1+r_1+\cdots +r_k)} \\ \\
 & & \displaystyle \times \Big(\sum_{j=1}^{k}s_j(n)\prod_{i\in Y\backslash\{j\}}(z-z_i)\Big)
\prod_{j\in X\backslash Y}(z-z_j)\, \, .
 \end{array}
\end{displaymath}
We can write
\begin{displaymath}
f_n(z) = \prod_{j=1}^{k}(z-\zeta(j,n)) \, \, , z \in \mathbb{C} \, ,
\end{displaymath}
where $\zeta(1,n),\cdots,\zeta(k,n)$ denote the zeros of $nP'q_n + Pq_n'$ such that $\zeta(j,n) \notin \{z_1,\cdots,z_k\}$ for $j = 1, \cdots,k$. Clearly $\zeta(j,n) \not= a$ for $j = 1,\cdots,k$. Let
\begin{displaymath}
f(z) = \prod_{j = 1}^{|X|}(z-\zeta_j)\prod_{j\in Y\backslash X}(z-z_j) \, \, ,
\end{displaymath}
for $z \in \mathbb{C}$.
By (iv) there is a  zero $\zeta_j$  of $P'$ such that $|a-\zeta_j| < 1$. Let $0 < \rho < 1 -|a-\zeta_j|$. Since by (iii) $f_n \to f$ uniformly in $\overline{D(0,1)}$, there exists $n_0$ such that
 for all $n \geq n_0$, there exists $j_n \in \{1,\cdots,k\}$ such that $\zeta(j_n,n) \in D(\zeta_j,\rho)$. Then we have
\begin{displaymath}
|a-\zeta(j_n,n)| \leq |a-\zeta_j| +|\zeta_j - \zeta(j_n,n)| < 1.
\end{displaymath}
The proof of the theorem is complete.

\bigskip

\noindent{\sl Remark 2.} In the first part of the proof the strong inequality in (iv) can be replaced by a weak inequality, together with the weak inequality in the statement. In such case, (iv) is just requiring that the Sendov conjecture is true for $P$ with respect to $a$.

\bigskip

\bigskip

\noindent{\sl Remark 3.} In the setting of Theorem 3 suppose that  $a = 1$ and $s_j(n) = 0$ for $j = 1,\cdots,k$. Then by Corollary 1 there exists $\zeta \in \overline{D(0,1)}$ such that $\zeta \not= a$, $p'(\zeta) = 0$ and $|\zeta - 1/2| \leq 1/2$.

\bigskip

\noindent{\sl Remark 4.} It should be noticed that (iv) is satisfied when $a < 1$ is closed to 1: See \cite{mi}. If $a = 1$ (iv) also holds unless $P(z) = z^{k+1}-1$ (see \cite{ru}). Moreover, unless $P(z) = z^{k+1}-1$ as before, (iv) is satisfied when $\sum_{j\in X}r_j$ is sufficiently large (see \cite{t} Remark 5.8). Clearly (iv) is also satisfied if $z_j$ is a multiple root and $|a-z_j| < 1$ for some $j \in X$.

\bigskip

\noindent{\sl Acknowledgment.} The author would like to thank the referee for the useful comments and suggestions.

\end{document}